\documentclass[11pt,a4paper,leqno]{amsart}
\usepackage{amsmath,amsthm,amsfonts,amscd,latexsym,amssymb}
\usepackage[all]{xy}
\newcommand{\ZZ}{\mathbb{Z}}

\newcommand{\CC}{\mathbb{C}}
\newcommand{\PP}{\mathbb{P}}
\newcommand{\NN}{\mathbb{N}}

\newcommand{\HH}{\mathbb{H}}

\newcommand{\QQ}{\mathbb{Q}}
\newcommand{\RR}{\mathbb{R}}

\newcommand{\Ker}{{\rm Ker}}

\newcommand{\Aut}{{\rm Aut}}

\newcommand{\Aff}{{\rm Aff}}

\newcommand{\Inn}{{\rm Inn}}
\newcommand{\Stab}{{\rm Stab}}
\newcommand{\Fix}{{\rm Fix}}
\newcommand{\Per}{{\rm Per}}
\newcommand{\SL}{{\rm SL}}
\newcommand{\GL}{{\rm GL}}
\newcommand{\PSL}{{\rm PSL}}

\newcommand{\free}[1]{\hat F(#1)}
\newcommand{\msC}{{\mathcal M}}
\newcommand{\mrC}{M}

\newcommand{\topo}{{\rm top}}
\newcommand{\sms}{\smallsetminus}

\newcommand{\pni}{{\PP^1_{0,1,\infty}}}

\DeclareMathOperator{\Sp}{Spec}

\newcommand{\ve}{{\varepsilon}}
\newcommand{\ol}{\overline}
\newcommand{\ul}{\underline}

\newcommand{\Centr}{{\rm Centr}}
\newcommand{\orb}{{\rm orb}}

\newtheorem{Defi}{Definition}[section]
\newtheorem{Rem}[Defi]{Remark}
\newtheorem{Lemma}[Defi]{Lemma}
\newtheorem{Cor}[Defi]{Corollary}
\newtheorem{Thm}[Defi]{Theorem}

\newtheorem{Prop}[Defi]{Proposition}
\usepackage[]{graphicx}
\begin{document}

\keywords{Teichmueller curve, Hurwitz space, faithful Galois action,
$\widehat{GT}$-relation}
\subjclass[msc2000]{14H30, 32G15}

\title{Teichm\"uller curves, Galois actions and 
$\widehat{GT}$-relations}


\author{Martin M\"oller}
\begin{abstract}
Teichm\"uller curves are geodesic discs in Teichm\"uller space
that project to algebraic curves $C$ in the moduli space $\mrC_g$. 
Some Teichm\"uller curves can be considered as components of Hurwitz
spaces. We show that the absolute Galois group
$G_\QQ$ acts faithfully on the set of these embedded curves.
\newline
We also compare the action of $G_\QQ$ on $\pi_1(C)$ with
the one on $\pi_1(\mrC_g)$ and obtain a relation 
in the Grothendieck-Teichm\"uller group, seemingly independent
of the known ones.
\end{abstract}
\maketitle                   

\section*{Introduction}

Consider a complex geodesic for the Teichm\"uller metric
$\tilde{\j}: \HH \to T_g$ from the upper
half plane to Teichm\"uller space. These geodesics are generated
by a pair $(X,q)$ of a Riemann surface $X$ of genus $g$
and a quadratic differential
$q$. The (rare) examples where the stabilizer in the mapping class
group of $\tilde{\j}$ is
a lattice $\Gamma \subset \Aut(\HH)$ are called {\em Teichm\"uller curves}.
\par
A particular case of these geodesics with lattice stabilizer
can be described as follows:
Take an assemblage of squares
of paper and glue them along their edges to a surface without
boundary, such that at each vertex abut an even number of squares.
If we provide the squares with a complex structure and glue
the local quadratic differentials $dz^2$, we obtain a pair $(X,q)$.
\newline
This description lead Lochak (\cite{Lo03}) to baptise them
{\em origamis}. If the glueing is just by translations (and
not by $(-1)$ composed by a translation) 
the origami is said to be {\em oriented}. 
Oriented origamis are also known as square-tiled coverings
or non-primitive Teichm\"uller curves. Lochak remarked that
the corresponding {\em origami curves} $j:C = \HH/\Gamma \to M_g$ in
the moduli space of curves are defined over
number fields and hence interesting not only in dynamical systems
but also from a number theoretical viewpoint.
\newline
We study two examples (section \ref{Exams}) of origamis:
the smallest example (called $L(2,2)$) where $X$ has genus greater than one 
and the smallest example (called ${\mathcal S}_2$) with a property relevant for
the $\widehat{GT}$-comparisons, see below. In both cases
we explicitly write down the equation of the origami curve 
in moduli space. This possibility seems rather unexpected from
the 'geodesic' viewpoint.
\par
The above definition suggests that -- in the oriented case -- origami 
curves are in fact (the analytic version of) some Hurwitz spaces for
coverings of elliptic curves, ramified at most over $\infty$. 
It is very natural to consider the geodesic discs then in
$\mrC_{g,[n]}$ with the $n$ preimages of $\infty$ marked but unordered.
The Hurwitz viewpoint 
reproves that origami curves and the map $j$ are defined over 
number fields. It implies that there is an action of the absolute Galois group
$G_\QQ$ on the set of origami curves. 
Our first main result is that this action is faithful in the following sense
(we abbreviate $\PP^1 \sms
\{0,1,\infty\}$
by $\PP^*$ throughout):
\par
{\bf Theorem \ref{faithfulThm}} {\it 
For each $\sigma \in G_\QQ$ there is an origami curve $C$ isomorphic
to $\PP^*$, such that $\sigma$
acts non-trivially on the map $j(C) \to \mrC_{g,[n]}$, 
more precisely such that $j(C) \neq  j^\sigma(C)$. }
\par
There are examples of origami curves, that do not have genus $0$. 
(\cite{Sm03}). It would be interesting to know, if 
the $G_\QQ$-action is faithful
on the (abstract, not embedded) origami curves.
\par
\medskip 
Let us give an overview over Grothendieck-Teichm\"uller theory.
The technical details will be explained in section \ref{onGT}. 
\newline
The absolute
Galois group $G_\QQ$ acts on the algebraic fundamental group
$\pi_1(\PP^\ast_{\ol{\QQ}})$, which is isomorphic to the profinite
free group $\widehat{F_2}$ in two generators. This action can be
described in terms of pairs $(\lambda_\sigma, f_\sigma) \in \widehat{\ZZ}
\times (\widehat{F_2})'$ with a suitable composition law. The
subgroup of pairs satisfying three equations (due to Drinfel'd)
is the Grothendieck-Teichm\"uller group $\widehat{GT}$. By
construction it contains $G_\QQ$ and the question is whether
it is strictly larger than $G_\QQ$ or not. One hence would like
to find a set of equations that singles out precisely $G_\QQ$.
\newline
Where should these relations come from? Consider a category
${\mathcal C}$ of varieties (or stacks) over $\QQ$ with 
some morphisms defined over $\QQ$. 
For a large enough ${\mathcal C}$ only $G_\QQ$ acts 
equivariantly on all $\pi_1({\mathcal C})$ by a result of F. Pop. 
\newline
But only for a few types of 'geometric' morphisms one is 
able to express the equivariance in terms of $(\lambda, f)$. 
Among these are some coverings of curves (\cite{NaSc00}, \cite{NaTs99}), 
and 'natural' morphisms between moduli spaces of curves 
(\cite{HaLoSc00}, \cite{NaSc00}, \cite{Sn02}). Maps from
curves to moduli spaces should be considered next.
\newline
The map of some origami curves to moduli space are 
defined over $\QQ$ and 'geometric' in this sense. We need that 
the curve $C \cong \PP^*$ and that it passes through a 
maximally degenerate point of the (compactified) moduli space. 
The smallest such origami is the ${\mathcal S}_2$.  
\newline
Let $\alpha_1,\ldots,\alpha_5$ denote
the standard generators (see section \ref{CompGalAct}) of the 
profinite mapping class group $\widehat{\Gamma_{2,0}}$. For $\sigma 
\in G_\QQ$ such that the
Kummer cocycle $\rho_2(\sigma)=0$ we obtain the following relation:
\par
{\bf Theorem \ref{newGTrel}'\ } {\it
The element $(\lambda,f) \in \widehat{GT}$ respects
the Galois actions on the morphism $j: C^\orb \to \msC_2$ induced from
the two-steps origami ${\mathcal S}_2$ (see figure $4$) 
if and only if
$$ 
f(\alpha_3, (\alpha_1^2 \alpha_2)^4)f(\alpha_1^2,\alpha_2^2)
f(\alpha_5^2,\alpha_4^2)
f(\alpha_2\alpha_4,\alpha_1^2\alpha_3 \alpha_5^2)  = 1
$$
holds in $\widehat{\Gamma_{2,0}}$. 
The elements $(\lambda,f)$  satisfying this relation form a 
subgroup of $\widehat{GT}$
containing $G_\QQ$.
}
\par
In section \ref{CompGalAct} we prove the refined relation for all $G_\QQ$.
As for all the relations in $\widehat{GT}$ recently discovered
it is not known, if the subgroup defined by this relation is
properly between $G_\QQ$ and $\widehat{GT}$. 
\par
The author thanks Pierre Lochak and Leila Schneps a lot for
introducing him to this subject and a lot of support. He also 
thanks W.~Herfort and H.~Nakamura for their suggestions.
\par

\section{Teichm\"uller curves}

A holomorphic quadratic differential $q \neq 0$ 
on a Riemann surface $X$ of genus $g$ determines on $X$ minus
the set of zeroes of $q$ an atlas of open charts, whose transition 
functions are of the form $z \mapsto \pm z + c$.
Such an atlas is called a {\em flat structure} on $X$. Conversely
a flat structure determines a quadratic differential by
glueing the local $dz^2$'s.
This correspondence is given in more details in \cite{Lo03} Ch.\ 2.
\newline15:58
There is a natural $\SL_2(\RR)$ action on pairs $(X,q)$ by 
postcomposing the charts of the flat structure with the 
linear map. The stabilizer of a pair $(X,q)$ is precisely
${\rm SO}_2(\RR)$. We define the {\em Teichm\"uller space} ${\mathcal T}_g$ 
as the space of Riemann surfaces plus an isotopy class of 
orientation-preserving diffeomorphism (a {\em Teichm\"uller marking}) 
to a reference surface $\Sigma_g$. 
If we choose a Teichm\"uller marking on $X$, the 
action $\SL_2(\RR) \cdot (X,q)$ yields 
a geodesic curve (by Teich\-m\"uller's theorems) 
$\tilde{\j}: \HH \to {\mathcal T}_g$.
We denote by $\msC_{g}$ the moduli stack of curves and
the corresponding coarse moduli space by non-calligraphic
letters, i.e.\ $\mrC_g$.
\par
Consider the projection of a geodesic $\HH \to {\mathcal T}_g$ 
to $\mrC_{g}$. 
We consider the moduli space in the first two sections
in the analytic category. More precisely we should write 
$(\mrC_{g})_{\CC}^{{\rm an}}$, etc. Let $\Gamma_{g,n}$ be the
mapping class group of Riemann surfaces of genus $g$ with $n$
punctures. We define $\Stab(\tilde{\j}) \subset \Gamma_{g,0}$
to be the (setwise) stabilizer of $\tilde{\j}(\HH)$ and 
$\Aut(\tilde{j})$ to be the pointwise stabilizer. We summarize
this by an exact sequence
$$ \xymatrix{
1 \ar[r] & \Aut(\tilde{\j}) \ar[r] & \Stab(\tilde{\j}) \ar[r] & 
\ol{\Stab(\tilde{\j})} \ar[r] & 1. }$$
The map $\tilde{\j}$ descends to a map 
$$j: C := \HH/\ol{\Stab(\tilde{\j})} \to \mrC_g.$$
\begin{Defi} This map $j$ is called a {\em Teichm\"uller curve} 
if\ \ $\ol{\Stab(\tilde{\j})}$
is a lattice in $\Aut(\HH)/\{\pm 1\}$.
\end{Defi}
\par
Later on it will sometimes be natural to fix say $n$ marked points, 
consider $\tilde{\j}: \HH \to {\mathcal T}_{g,n}$ and consider
the stabilizer in $\Gamma_{g,n}$ (or in $\Gamma_{g,[n]}$, if
we allow permutation of the marked points). We then call the
groups $\Stab(\tilde{\j},n)$ (or $\Stab(\tilde{\j},[n])$) etc.
\par
We can also describe this group using the flat structure
defined by $(X,q)$. Denote by $\Aff^+(X,q)$ the group of
orientation preserving affine diffeomorphisms of $X$, i.e.\
diffeomorphisms 
which are affine with respect to the charts of the flat
structure determined by $q$. Associating with $\varphi
\in \Aff^+(X,q)$ the matrix part of the affine maps yields
a well-defined map $D$ to $\PSL_2(\RR)$. We denote the image of $D$ 
by $\PSL(X,q)$. We obtain an exact sequence
$$ \xymatrix{
1 \ar[r] & \Aut(X,q) \ar[r] & \Aff^+(X,q) \ar[r]^D & 
PSL(X,q) \ar[r] & 1,}$$
where $\Aut(X,q)$ are the (conformal) automorphisms of $X$ preserving $q$.
\par
We will also consider the subgroup
of $\Aff^+(X,q)$ that fixes $n$ points (resp.\ up to
permutation). We will denote these groups by 
$\Aff^+(X,q,n)$ (resp.\ $\Aff^+(X,q,[n])$) and their images under
$D$ by $\PSL(X,q,n)$ (resp.\ $\PSL(X,q,[n]$)). Affine diffeomorphisms
in $\Aff^+(X,q,[n])$ are called {\em balanced}.
\par
\begin{Rem} {\rm i) 
The groups $\ol{\Stab(\tilde{\j})}$ and $\PSL(X,q)$ are closely related,
namely (see \cite{Mc02} Prop.\ 3.2) if we identify $\Aut(\HH)/\{\pm 1\}$
with $\PSL_2(\RR)$, we have
$$ \ol{\Stab(\tilde{\j})} =  R \cdot  \PSL(X,q) \cdot R, \quad 
\text{where} \quad R = \left(\begin{array}{cc} -1 & 0 \\ 0 & 1 \\ \end{array}
\right).$$  The same holds for the
corresponding groups with marked points. \newline
ii) While 
$\ol{\Stab(\tilde{\j})}$ is best suited for moduli problems, we
will use the definition of $\PSL(X,q)$ for calculations of
these groups, because affine diffeomorphisms are easily visualized.}
\end{Rem}
\par

\section{Origamis} \label{Origamidef}

We consider now a special case of Teichm\"uller curves:
\par
\begin{Defi} An origami is a finite set of squares (say unit
squares in $\RR^2$) glued together along their edges to a surface without
boundary, such that at each vertex abut an even number of squares.
\end{Defi}
\par
Using also open charts covering the glueing edges and identifying
$\RR^2$ with $\CC$, this construction defines a Riemann surface $X$
with a flat structure. The local $dz^2$'s glue
together to a quadratic differential $q$, which is holomorphic
except for simple poles at vertices where only $2$ squares abut.
\par
We can distinguish two cases: If $q = \omega^2$ is a square of
a one-form $\omega \in H^0(X,\Omega_X)$ the origami is said to 
be {\em oriented}, otherwise {\em non-oriented}.
\newline 
If the origami is oriented, the number of squares abutting at 
each vertex is divisible by $4$. This condition is not
sufficient, as shown by the following example, where $A$
is glued to $A'$ etc.\ along the orientation of the arrow:
\begin{figure}[h] \label{nonor}
\centerline{\includegraphics{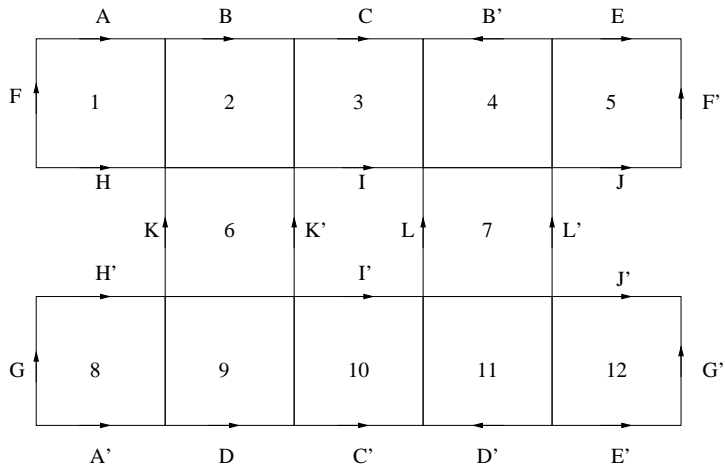}}
\caption{A non-oriented origami without poles}
\end{figure}
\par
Note however, that if the origami is non-oriented, there is a canonical double
covering (ramified precisely over the points, where $2$ squares abutt), 
which is an oriented origami. 
\newline
{\em In the sequel we will treat only the oriented case}. A neccessary
and sufficient condition for an origami to be oriented is that
the transition functions between the charts consist only
of translations. In origami language this means that upper (resp.\ left)
edges should be glued to lower (resp.\ right) edges preserving
global orientation.
\par
\begin{Defi}
An origami in the above sense defines an 
unramified covering $\pi:X^* \to E^*$ of a torus punctured
at $\infty$ and a ramified covering $X \to E$, also denoted by $\pi$.
We call this an {\em origami covering}.
\end{Defi}
\par
Note that by identifying $\RR^2$ with $\CC$ we used the elliptic
curve with $j(E) = 1728$ in the first definition. This choice plays
no role in the sequel, because the geodesic curve generated by $(X,q)$
is independent of this choice.
\par
The fact that origamis indeed define Teichm\"uller curves, follows
from the result of \cite{GuJu00}, that $\ol{\Stab(\tilde{\j})}$
is commensurable with $\PSL_2(\ZZ)$. This will also follow from
the Hurwitz space description in section \ref{OrigIsHur}.
\par
\begin{Rem} \label{MonoRem}
{\rm We will specify the (unramified) covering
$\pi$ of degree $d$ by its monodromy: Fix two generators
$a$ and $b$ of $\pi_1(E^*)$ and their images under
the monodromy map $m:\pi_1(E^*) \to S_d$. These images
determine $\pi$. Note that simultaneous conjugation in $S_d$
(i.e.\ renumbering the preimages of a basepoint) gives
the same covering. Note also that topologically different
coverings may lead to the same origami curve.
}\end{Rem}
\par
The quadratic differential $q$ on $X$ is obtained using
the origami covering as $q = \omega^2$, where $\omega = \pi^* \omega_E$
and $\omega_E$ is the unique (up to scalar multiple)
holomorphic differential on $E$. 
\par
In the case of origamis it is natural to consider the
groups $\Stab(\tilde{\j},[n])$ (or $\Stab(\tilde{\j},n)$), where the marked
$n$ points are the preimages $\pi^{-1}(\infty)$. This is for two reasons:
First, this finite set of additional marked points replaces 
$\Stab(\tilde{\j})$ by a subgroup of finite index (see \cite{GuJu00}).
Note that not any set of additional marked points has this property, 
see \cite{GuHuSc03}.
Second, the following definitions are best suited with the
Hurwitz space interpretation in section \ref{OrigIsHur} and
yield that $\ol{\Stab(\tilde{\j})}$ is contained in $\PSL_2(\ZZ)$.
\par 
\begin{Defi}
The group $\Gamma(\pi):= \ol{\Stab(\tilde{j},[n])}$ is called the
{\em affine group} of the origami. The quotient $C(\pi) = \HH/\Gamma(\pi)$ 
(or simply $C$) is called the {\em origami curve}. If we consider it
as orbifold quotient $C^{\rm orb}(\pi) := \HH/\Stab(\tilde{j},[n])$ it
is called the {\em orbifold origami curve}.
\end{Defi}
\par
The map $\HH/\PSL(X,q) \to M_g$ is generically injective
and actually injective up to finitely many 
normal crossings, because it is the image of a geodesic locus under the
quotient map of a discrete group. But the map
$C \to M_g$ might be (see the following example) a composition
of a covering and a generically injective map.
If we want to reestablish injectivity, we will consider the 
origami curve in $M_{g,[n]}$, the moduli space of curves
with $n$ non-ordered points. 
\newline
To give an example that
$\Gamma(\pi)$ is a proper subgroup of $\ol{\Stab(\tilde{\j})}$ 
consider the following thick $L$ (see figure $2$:
The sides $L_i$ are glued with $R_i$ and $D_i$ with $U_i$, 
preserving the global orientation. One checks that for this origami the matrix
$\left( \begin{array}{cc}
0 & -1/2 \\
2 & 0 \\
\end{array}\right)$ is in $\PSL(X,q)$ but neither in $\Gamma(\pi)$
nor in $\PSL_2(\ZZ)$.
\par
\medskip
\begin{figure}[h] \label{figlmn}
\centerline{\includegraphics{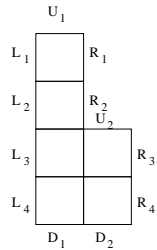}}
\caption{The thick $L$}
\end{figure}
\par
\medskip
\begin{Rem} {\rm
The affine group also admits the following description (see
also \cite{Sm03}):
\newline
One can always lift a balanced affine diffeomorphism
$\varphi \in \Aff^+(X,q,[n])$ to its universal cover $\HH$.
Indeed one can lift an affine diffeomorphism locally to an
unramified cover and non-trivial paths provide the
obstruction to do this globally. Denote by $q_\HH$
the quadratic differential on $\HH$ obtained by pullback
of the square of $\omega_E$ 
via the universal covering map $\pi_\infty: \HH \to E$.
One has a natural morphism
$$ \ast: \left\{ \begin{array}{lcl}
\Aff^+(\HH,q_\HH) &\to& \Aut^+(\pi_1(E^*)) \\
\varphi & \mapsto & \varphi_*:=(f \mapsto \varphi^{-1} \circ f \circ
\varphi)
\end{array} \right.,$$
where we consider $\pi_1(E^*)$ as $\Aut(\HH/E)$ and composition
is meant in $\Aut(\HH)$. The 'plus' of $\Aut^+(\pi_1(E^*))$
denotes the preimage of $\SL_2(\ZZ)$ under the quotient map
by $D:\Aut(\pi_1(E^*))/ \Inn(\pi_1(E^*)) \to \GL_2(\ZZ)$. 
Denote by $\Aut_i(\HH/E)$ the automorphisms of $\HH$ over
the identity or the elliptic involution of $i$. 
The right and the left vertical morphism of the
commutative diagram with exact rows 
$$ \xymatrix{
1 \ar[r]& \Aut_i(\HH/E)  \ar[r]
 \ar[d] & \Aff^+(\HH,q_\HH) \ar[r] \ar[d]^\ast & \PSL_2(\ZZ) \ar[r]
\ar[d] & 1 \\
1 \ar[r] & \Inn(\pi_1(E^*)) \rtimes \langle i_* \rangle
\ar[r]& \Aut^+(\pi_1(E^*)) \ar[r]
& \PSL_2(\ZZ) \ar[r] & 1 \\
}
$$
are isomorphisms, hence $\ast$ is an isomorphism, too.}
\end{Rem}
\par
In this language the affine group $\Gamma(\pi)$ is the
image under $D$ of the subgroup of $\Aut^+(\pi_1(E^*))$ that
fixes $\pi_1(X)$ (as set).
\par
With this description the following Lemma is obvious:
\par
\begin{Lemma} \label{Normalintermed}
Suppose an origami covering factors as $\pi =\psi \circ
\tilde{\pi}: X \to Y \to E$. Denote the corresponding unramified
coverings by $X^* \to Y^* \to E^*$. 
If the fundamental group $\pi_1(Y^*)$ is characteristic
in $\pi_1(E^*)$ any balanced affine diffeomorphism $\varphi: X \to X$ 
descends to a balanced affine diffeomorphism $\varphi_Y: Y \to Y$.
\end{Lemma}
\par

\begin{Rem} \label{Strebelrem}
{\rm Recall that a holomorphic quadratic differential $q$
is called {\em Strebel}, if its horizontal trajectories
are compact or connect two zeroes of $q$. A direction $e^{2\pi i
\theta}$ for $\theta \in \RR$ is called Strebel for $q$, if $e^{2\pi i
\theta}q$ is Strebel. A Strebel differential decomposes the
Riemann surface into cylinders swept out by trajectories. With our
convention for origamis ($\RR^2 \cong \CC$ and $q$ is made by
local $dz^2$), the direction $e^{2\pi i\theta}$ is Strebel
if and only if $\theta \in \QQ$. We will tacitly assume this
in the sequel.}
\end{Rem}
\par
\begin{Rem} \label{Monocylrem}{\rm
If the origami is given by its monodromy $m(a), m(b) \in S_d$ 
the horizontal trajectories decompose $X$ into $c_a$ disjoint unit height
cylinders, where $c_a$ is the number of cycles of the permutation $m(a)$.
These cylinders lie in $c_a^{{\rm max}} \leq c_a$ maximal cylinders.
Denote by $\alpha_i$, $i=1,\ldots, c_a^{{\rm max}}$ the core curves of
the maximal cylinders. If both sides of a unit height cylinder contains
a zero of $\omega$, it is of course a maximal cylinder.
\newline
Consider the family of curves 
$$ X_t  := \left( \begin{array}{ll} e^t & 0 \\ 0 & e^{-t} \\ \end{array}
\right) \cdot X.$$ 
By Thm.~3 in \cite{Ma75} precisely the hyperbolic lengths
of the homotopy classes of $\alpha_i$ tend to zero (as $t$ tends
to infinity). Hence
the family of smooth curves $X_t$ tends for $t \to \infty$ 
to the stable curve obtained by 'pinching' $\alpha_i$ to
nodes.}
\end{Rem}
\par

\section{Origami curves are components of Hurwitz spaces}\label{OrigIsHur}

Origami  curves were defined in the previous sections as quotients
$C = \HH/\Gamma(\pi)$ of discs in Teichm\"uller space, hence in the
analytic category. But the differential $\omega$ (or $q=\omega^2$)
was defined by covering data of an elliptic curve (e.g.\ $E_{1728}$).
One expects that the (orbifold) origami curve in the analytification
of a component of an algebraic Hurwitz (stack resp.) space. We prove here
that this is indeed the case.
\par
We start with generalities on (algebraic) Hurwitz stacks. Roughly speaking
they parametrize isomorphism classes of coverings. 
Isomorphism here means isomorphism over a {\em fixed} 
base curve. Notations 
in this section follow Wewers (\cite{We98}). We only need
covers of smooth schemes and we assume all schemes to be 
schemes over $\QQ$. 
\par
Let $f:{\mathcal E}\to{\mathcal M}$ denote the universal family of
curves over the smooth stack ${\mathcal M}$. Fix a ("ramification") 
divisor ${\mathcal D}/{\mathcal M}$. Of course we have elliptic
curves, maybe with additional structures in mind.
\newline
For any ${\mathcal M}$-scheme $S$ let 
${\mathcal H}_{\mathcal E}(S)$
be the category of finite covers $X \to E_S$ of fixed
degree $d$ and fixed genus $g$ of $X$, ramified over 
$D_S := {\mathcal D} \times_{\mathcal M} S$, where 
$E_S := {\mathcal E} \times_{\mathcal M} S$. 
(Maybe the notation ${\mathcal H}_{{\mathcal E},g,d}$ would be more
precise.) A morphism from
$X \to E_S \to S$ to $X' \to E_{S'} \to S'$  
is a morphism $S \to S'$ plus a morphims $X \to X'$ over
the induced morphism $E_S \to E_{S'}$. We cite Th.\ 4.1.2 from \cite{We98}:
\par
\begin{Prop} \label{Wewprop}
The Hurwitz stack ${\mathcal H}_{\mathcal E}$ is a smooth stack over 
$\QQ$, {\'e}tale over ${\mathcal M}$. The stack
${\mathcal H}_{\mathcal E}$ has a coarse moduli space $H_{\mathcal E}$, which is
also defined over $\QQ$.
\end{Prop}
\par
If we replace in the above definition ${\mathcal M}$ by
${\mathcal M}_\CC^{\rm an}$ we obtain an analytic stack 
${\mathcal H}_{\mathcal E}^{\rm an}$, 
whose coarse moduli space is $(H_{\mathcal E})_\CC^{\rm an}$. In the
sequel we do not distinguish between ${\mathcal E}$ and
${\mathcal E}_\CC$ for simplicity of notation.
\par
\begin{Prop} \label{HurProp}
Let ${\mathcal E}$ be the universal family over
the moduli stack ${\mathcal M}_{1,1}$ and ${\mathcal D}$ the divisor corresponding
to the section. Let ${\mathcal H_\pi}$ be the connected component of 
${\mathcal H}_{\mathcal E}$ which contains
$\pi: X \to E$.
\newline
Then the orbifold origami curve
$C^{\rm orb} = C^{\rm orb}(\pi)$ coincides
(as functor (analytic spaces) $\mapsto$ (sets))
with $({\mathcal H_\pi})_\CC^{\rm an}$.
\end{Prop}
\par
{\bf Proof:} 
Consider the analytic stack
${\mathcal TH}_{\mathcal E}$ whose objects (over $S$)
are coverings $X \to E_S \to S$ in $({\mathcal H}_{\mathcal E})_\CC^{\rm an}(S)$ 
with compatible Teichm\"uller markings on $X$ and $E_S$.
Because  ${\mathcal T}_g$ is \'etale over $(\msC_g)^{\rm an}_\CC$
and because of Prop.\ \ref{Wewprop}
the functor 'forget $X$' exhibits 
$({\mathcal TH}_{\mathcal E})$ as an \'etale
cover of ${\mathcal T}_{1,1} \cong \HH$. Hence it consists
of several disjoint copies of $\HH$, one of which
contains $\pi: X\to E$. We call this connected component
${\mathcal TH}_\pi$ and we let $({\mathcal H}_\pi)^{\rm an}_\CC$ be 
the corresponding component of $({\mathcal H}_{\mathcal E})_\CC^{\rm an}$ 
(after forgetting the Teichm\"uller marking). $({\mathcal H}_\pi)^{\rm an}_\CC$ 
stems from a uniquely determined  connected component ${\mathcal H}_\pi$
of ${\mathcal H}_{\mathcal E}$.
\par
It remains to check that $({\mathcal H}_\pi)^{\rm an}_\CC$ is
the quotient stack ${\mathcal TH}_\pi/ \Stab(\tilde{\j},[n])$.
This easily follows from 
unwinding definitions, noting that Deck
transformations respect the set $\pi^{-1}(\infty)$.
\hfill $\Box$
\par
\begin{Cor}
Origami curves are geometric components of a Hurwitz space
defined over $\QQ$. Hence they are defined over number
fields and there is a natural $G_\QQ$-action on the
set of origami curves.
\newline
The morphism $j:C \to M_g$ (and its orbifold version)
is defined over a number field and there is a natural 
$G_\QQ$-action on the set of embedded origami curves.   
\end{Cor}
\par
{\bf Proof:} For the second claim note that the
forgetful functor ${\mathcal H}_{\mathcal E} \to {\msC_g}$ is
defined over $\QQ$. Hence the morphism between the
geometric components of the (coarsely) representing 
schemes are defined over finite extensions of $\QQ$.
\hfill $\Box$.
\par
\begin{Rem} {\rm There are obviously substacks of
${\mathcal H}_{\mathcal E}$ which are defined over $\QQ$,
for example the Hurwitz spaces with fixed monodromy
(see e.g.\ \cite{We98}). Each Galois invariant
additional structure describes a substack defined
over $\QQ$. Finding components of ${\mathcal H}_{\PP^*}$
that are irreducible over $\QQ$ has been studied
in the context of dessins d'enfants under the
name of giving a {\em complete list of Galois
invariants}. One can ask the same question for ${\mathcal H}_{\mathcal E}$.
\newline
Here is a list (certainly not complete) of
Galois invariants known to the author. The given
references consider the invariants from quite different viewpoints. 
\newline
$\bullet$ Monodromy groups and ramification indices, or 
more generally the Nielsen classes are Galois invariant.
\newline
$\bullet$ If the ramification indices of $\pi$ over $\infty$
are all odd, the parity of the spin structure (see
\cite{KoZo02}) is a Galois invariant.
\newline
$\bullet$ If $g(X)=2$ there is a $2$-division
point $\mu_2$ in $E$ such that $\pi+\mu_2: X \to E$
is equivariant with respect to the (hyper)elliptic
involutions on $X$ and $E$ (see \cite{Ka03}). 
The property whether or not $\mu_2 = 0$ is a Galois invariant.
}
\end{Rem}
\par

\section{Two examples} \label{Exams}

In this section we examine two origamis, the smallest
origami which is not an elliptic curve and the smallest
one that can be used for $\widehat{GT}$-considerations.
In both cases we explicitely describe the equation 
of the geodesic curve.
\par
\subsection{The $L(2,2)$}

We now study the simplest origami which is not an elliptic
curve, i.e. such that $\pi$ is not an isogeny. It is 
called $L(2,2)$ in \cite{Lo03}, see the left half of
figure $3$. The sides are glued 'naturally', i.e.\ $L_i$
with $R_i$ and $U_i$ with $D_i$. The affine group 
of this origami $\Gamma(L(2,2))$
contains the horizontal and vertical translation by $2$, 
hence the modular group $\Gamma(2)$ and also $$S = \left( 
\begin{array}{cc} 0 & 1 \\ -1 & 0 \\ \end{array} \right),
\quad \text{but not} \quad T = \left( \begin{array}{cc} 1 & 1 \\ 0 & 1 
\\ \end{array} \right).
$$
In fact $S$ corresponds to rotation by $90^\circ$. And if
$T$ were in $\Gamma(L(2,2))$ one could map $R_1$ to one
of the other vertical edges and extend this map to a
diffeomorphism, that locally looks like $T$. But this leads
to a contradiction in each of the cases. See \cite{Sm03}
for an algorithm to determine the affine group of an origami.
We conclude that $C \to \mrC_{1,1}$ is a degree $3$ cover. 
\par
We can illustrate Remark \ref{MonoRem} here: Fix $E$ and $a,b \in \pi_1(E)$.
The coverings with mono\-dromy $m(a)= (1)(23)$, $m(b)= (12)(3)$
(the $L(2,2)$) and with monodromy $m'(a)=(123)$, $m'(b)=(12)(3)$
are topologically different (i.e.\ one cannot obtain one
from the other by Deck transformations and renumbering of the
preimages). Nevertheless the origami curves coincide.
\par
We can see geometrically, that the corresponding origami
covering $\pi:= \pi_{L(2,2)}$ commutes with the hyperelliptic
involutions $h_X$ and $h_E$ on $X$ and $E$ respectively
(see figure $3$).
\par
\medskip
\begin{figure}[h] \label{figl22}
\centerline{\includegraphics{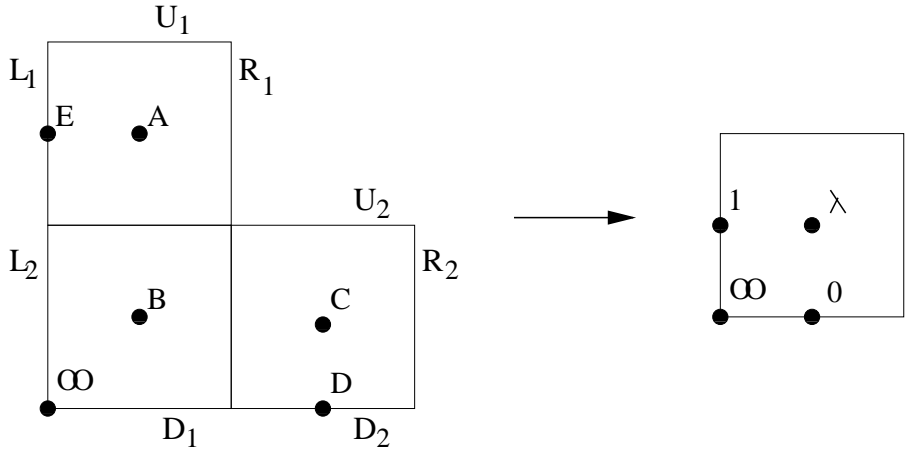}}
\caption{Weierstra\ss \ points of the $L(2,2)$}
\end{figure}
\par
When giving the equation of this family, we work for simplicity
over $\PP^\ast = \PP^1 \sms \{0,1,\infty \}$ (with coordinate $t$)
instead of its quotient by $S$ to have all Weierstra\ss\ points available.
\par
\begin{Prop} \label{eq22prop}
The equation of the family ${\mathcal X}/\PP^*$ is given by
$$y^2 = x(x-4)(P(x)-t),\, \text{where} \, P(x)= \frac{1}{4}x(x-3)^2.$$
together with the morphism $\pi(x,y) = (P(x), yQ(x))$, where
$Q(x)=(x-3)(x-1)/4$.
\end{Prop}
\par
{\bf Proof:}
Denote by $q:X \to X/h_X$ and $q_E: E \to E/h_E$ the quotient maps.
To determine the equations it is sufficient to find a map 
$P: X/h_X \to E/h_E$ of degree $3$ with the correct ramification behaviour: 
As the preimage of $\infty$ under $q_E \circ \pi$ is just one
point (also denoted by $\infty$ in the above figure), we may
suppose that $P$ is a polynomial, and hence $q_X(\infty) = \infty$.
Furthermore the preimage of $0$ (resp.\ $1$) under $q_E \circ \pi$
consists of a Weierstra\ss\ points of $X$ (namely $D$ resp.\ $E$) 
and two points
exchanged by $h_X$. Hence among $P^{-1}(0)$ (resp.\  among $P^{-1}(1)$)
there must be precisely one ramification point.
The polynomial $P(X)$ does the job and the images in $X/h_X$ of 
the Weierstra\ss\ points are $\infty$, the ramification points of $P$ over
$0$ and $1$ and the three preimages of some $\lambda \not\in \{0,1,\infty\}$.
\hfill $\Box$
\par

\subsection{The two steps} \label{Ex2Steps}

We will study the following origami, let's call it ${\mathcal S}_2$,
given by the 
permutations $m(a)=(12)(34)$ and
$m(b)=(1)(23)(4)$ or graphically by 'two steps'
\par
\begin{figure}[h] \label{figorigs2}
\centerline{\includegraphics{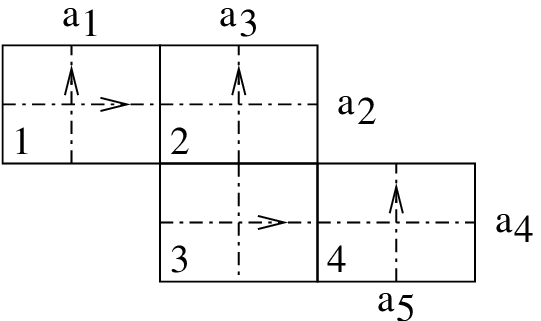} \ \ \ \ \ \ 
\includegraphics{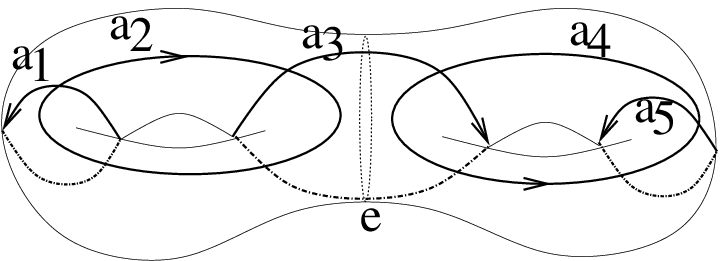}}
\caption{The origami ${\mathcal S}_2$}
\end{figure}
\par \noindent
with the natural identifications, i.e.\ glueing left sides to right
sides in the same row in the same orientation and glueing top
sides to bottom sides in the same column. (The surface on the right
does obviously not contain the origami grid. The loop $e$ is
added for later use.)
\par
The origami ${\mathcal S}_2$ is of genus $g=2$
and it is the origami of lowest degree, 
that has a Strebel direction with $3g-3$ cylinders (namely e.g.\
the vertical ones). This means (see Rem.~\ref{Monocylrem}) that the
origami curve passes through a maximally degenerate point
of (the compactification of)
$\mrC_2$ and will be important in section \ref{CompGalAct}.
Here {\em maximally degenerate} means that the stable curve
consists of a graph of $\PP^1$'s with three marked points
or normal crossings.
\par
For each elliptic curve $E$ the 
covering $\pi: X \to E$ is ramified at two points of order $2$ (but
not Galois). It factors as $\pi = \iota \circ \pi_1$ in a degree
$2$ covering $\pi_1 : X \to E_1$ and an isogeny $\iota: E_1 \to E$
of degree $2$. 
\par
\begin{Prop} \label{eqprop}
$C = C({\mathcal S}_2)$ is the modular curve 
$\mrC^{[2]}_{1,1} = \HH/\Gamma(2) =: \PP^*$.
$C$ parametrizes the curves of genus $2$ given by
$$ y^2 = ([4x(x-1)]^2 - (1-t))(4x^2 - 4x -1) = (P^2 - (1-t))(P-1),$$
where $t$ is a coordinate on $\PP^*$ and where $P(x)= 4x(x-1)$.
$t$ is normalised such that $t=0$ gives the maximally degenerate point.
\newline
The family of intermediate covers ${\mathcal E}_1 / \PP^*$  
(containing $E_1$) is given by 
$$ y^2 = (x-1)(x+1)(x^2 - (1-t)).$$
The morphisms between these curves are
$$\pi_1: (x,y) \mapsto (P(x), 2y(x-\frac{1}{2})))$$
and
$$\iota: (x,y) \mapsto (x^2, yx),$$
where the family of base curves is given by $ y^2 = x(x-1)(x-(1-t))$.
In particular $C$ and $j:C \to \mrC_2$ are defined over $\QQ$.
\end{Prop}
\par
{\bf Proof:} One notices from figure $4$ that $\Gamma(2)$ 
is contained in the affine group $\Gamma(\pi)$ 
of ${\mathcal S}_2$ and by inspection (using e.g.\ \cite{Sm03}) 
one finds, that it is not bigger. 
\newline
To determine the equations note that the hyperelliptic involution 
$h_X$ of $X$ is compatible with $\pi_1, \iota$
and the elliptic involution $h_{E_1}$ and $h_{E}$.
The Weiersta\ss \ points of these curves are as follows:
\par
\begin{figure}[h] \label{figwpons2}
\centerline{\includegraphics{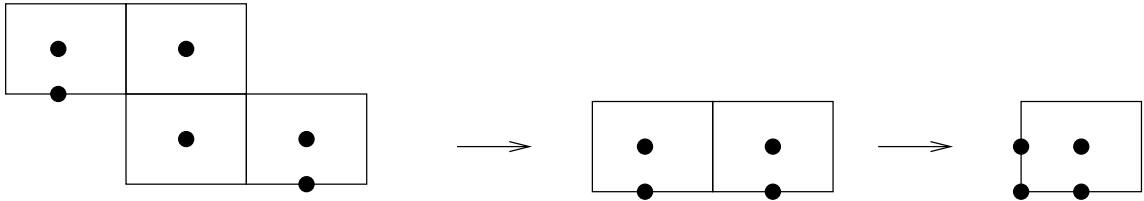}}
\caption{Weierstra\ss \ points}
\end{figure}
\par
To determine $\pi$ it is sufficient to find maps 
$X/h_X \to E_1/h_{E_1} \to E/h_E$ with the correct behaviour and
in fact the polynomials $P= 4X(X-1)$ and $X^2$ do the job.
\hfill $\Box$
\par
\begin{Rem} \label{orbs2}
{\rm The orbifold structure of $C^{\orb}({\mathcal S}_2)$
consists precisely of a globally acting group
$(\ZZ/2\ZZ)^2$ given by the hyperelliptic involution and
the automorphism $(x,y) \mapsto (1-x,y)$.
}\end{Rem}
\par

\section{The $G_\QQ$-action on oriented origamis is faithful}
\label{GalAct}
Fix $\sigma \in G_\QQ$ and let $K={\rm Fix}(\sigma)$ its fixed
field. We first prove faithfulness in a weak sense:
\par
\begin{Lemma} \label{weakfaithful}
For each $\sigma \in G_\QQ$ and each elliptic curve $E$ defined
over $K$ there is a covering $\pi: X \to E$, unramified
except over one point, such that 
$^\sigma \pi$ is not $\ol{\QQ}$-isomorphic to $\pi$.
\end{Lemma}
\par
{\bf Proof:} We will derive this from a corresponding
result concerning dessins d'enfants:
\newline
As shown in \cite{Sn94} Th.\ II.4 there exists a Belyi morphism 
$\beta: \PP_{\ol{\QQ}}^1 \to \PP_{\ol{\QQ}}^1$, i.e.\ unramified outside
$\{0,1,\infty\}$, such that $^\sigma \beta \not\cong \beta$.
We may also suppose that $\beta$ is pure, i.e. that
precisely the preimage of $1$ consists only of points ramified 
of order $2$.
\newline
Let $h: E \to \PP^1$ be the double cover ramified over 
$\{0,1,\lambda,\infty\}$.
We take the pullback of $\beta$ by the morphism $h$ and
call the desingularisation $\widetilde{\pi}:X \to E$. The morphism
$\widetilde{\pi}$ is ramified at most over $\{0,1,\lambda,\infty\}$
and hence after multiplication by $2$ the map $\pi = 
[2] \circ \widetilde{\pi}$ 
defines an origami.
\newline
Suppose there is a $\ol{\QQ}$-isomorphism
$\varphi: X \to \,^\sigma\! X$ with $\pi = \,^\sigma \pi \circ
\varphi$. We claim that this implies $\widetilde{\pi} = \,^\sigma 
\widetilde{\pi} \circ
\varphi$: As $E$ and $[2]$ are defined over $K$ the difference
morphism $\widetilde{\pi} - \,^\sigma \widetilde{\pi} \circ \varphi$ 
maps to the
finite kernel of $[2]$. We can conclude because $\beta^{-1}(1)$ and hence
$(\widetilde{\pi})^{-1}(h^{-1}(1))$ is 
distinguished by the purity hypothesis.
\newline
Let $\tau$ be the involution on $X$ induced by the elliptic
involution of $E$. Once we have shown that $\varphi \circ \tau 
=\, ^\sigma\! \tau \circ \varphi$, the isomorphism $\varphi$ 
descends to $\PP^1$ and gives the desired contradiction.
But both $\tau$ and $ \varphi^{-1} \circ\, ^\sigma\! \tau \circ \varphi$
are involutions with the same number of fixed points. Since $X/\langle \tau
\rangle$ the involutions have to coincide by the uniqueness of the 
hyperelliptic involution.
\hfill $\Box$
\par
Looking closer at the proof of Th.\ II.4 in \cite{Sn94}, 
we can choose the morphism $\beta$ in the above proof
such that its ramification over $1$ consists only of points
of order $2$, that $\beta$ is totally ramified over $\infty$ and
that the ramification behaviour over $0$ is different from these
two. 
Let $\msC_{1,1}^{[2]}$ denote the moduli stack of 
elliptic curves with level $2$ structure. We will use
the non-calligraphic letters (e.g.\ $\mrC_{1,1,}^{[2]} \cong \PP^\ast$)
for the corresponding coarse moduli spaces.
\par
\begin{Lemma} The orbifold origami curve  $C^\orb(\pi)$ as constructed in 
the preceding lemma starting with a Belyi morphism $\beta$, whose ramification
behaviour over $0$, $1$ and $\infty$ is pairwise distinct, 
is isomorphic to $\msC_{1,1}^{[2]}$, hence $C(\pi) \cong 
\PP^*$.
\end{Lemma}
\par
{\bf Proof:}
We continue with the construction and the notations from the proof 
of the above lemma. First we prove that $\msC_{1,1}^{[2]}$ 
surjects onto $C^\orb(\pi)$:
\newline
For this purpose we need to construct a family over 
$\msC_{1,1}^{[2]}$ of coverings
of the universal family ${\mathcal E} \to \msC_{1,1}^{[2]}$.
Let $h: {\mathcal E} \to \PP^1$ be the quotient by the elliptic
involution $h_E$ followed by the canonical projection $\PP^1 \times 
\msC_{1,1}^{[2]}
\to \PP^1$. (The quotient ${\mathcal E}/h_E$ is indeed a trivial
bundle as it contains $4$ disjoint (Weierstra\ss ) sections).
Take the desingularisation of $\PP^1 \times_{\PP^1}
{\mathcal E} \to \msC_{1,1}^{[2]}$, where the morphisms in the fibred
product are $\beta$ and $h$ respectively.
The singular locus of the fibre product is \'etale over the base
and thus the desingularisation is still a flat family of curves,
which we denote by ${\mathcal X} \to
\msC_{1,1}^{[2]}$.
Let $\pi: {\mathcal X} \to {\mathcal E}$ denote the composition of the
second projection with the multiplication by $2$.
We thus constructed a topologically locally trivial family of
coverings over $\msC_{1,1}^{[2]}$. 
This family contains by construction the covering 
$X \to E$ of the preceding lemma (which we also denoted by $\pi$),
hence does the job by Prop.\ \ref{HurProp}.
\par
It remains to exclude that $\msC_{1,1}^{[2]} \to C^{\orb}$
is a cover of degree greater than one or equiva\-lently
that $\PP^* \to C(\pi)$ is a cover of degree greater than one.
We have to exclude that the affine
group of $\pi$ is bigger than $\Gamma(2)$. 
\newline
Fix a fibre $X \to E^{[2]} \to E$ of $\pi: {\mathcal X} \to {\mathcal E}$
and call this covering also $\pi = [2] \circ \tilde{\pi}$. 
The subgroup of $\pi_1(E^\ast) = \langle x,y \rangle$
corresponding to $[2]$ is generated by $x^2, y^2, xy^2x$ and
$yx^2y$. It is obviously
characteristic and hence by Lemma \ref{Normalintermed} an 
affine diffeomorphism $\varphi: X \to X$ 
over $\ol{\varphi}: E \to E$ descends to an affine diffeomorphism
$\varphi^{[2]}:E^{[2]} \to E^{[2]}$. If $D(\ol{\varphi})$ is
not in $\Gamma(2)$, the morphism $\varphi^{[2]}$ has to
permute the $2$-division points of $E^{[2]}$. 
\newline 
But the fibres of $\tilde{\pi}$ over the $2$-division points 
are different by our hypothesis on $\beta$:
One fibre consists of unramified points that are all fixed by the
hyperelliptic involution of $X$. One consists of unramified points 
that are pairwise interchanged by the
hyperelliptic involution. One consists of $2$ points and
the last one has a ramification behaviour different from the above.
\newline
As an affine diffeomorphism has to preserve the ramification order
and fixed points of the hyperelliptic involution, this
leads to a contradiction.
\hfill $\Box$
\par
As usual let $d=\deg(\pi)$. We need one more topological lemma.
For the notation compare with the Remarks \ref{Strebelrem} and 
\ref{Monocylrem}.
\par
\begin{Lemma} \label{Belyiadjust}
We may suppose that the differential 
$\omega_X = \pi^* \omega_E = \tilde{\pi}^* \omega_{E^{[2]}}$
has in the horizontal (resp.\ vertical, diagonal) direction
$1$ (resp.\ $d/8$, resp.\ $r \not\in \{1,d/8\}$) maximal cylinders.
\end{Lemma}
\par
{\bf Proof:}
Consider $h: E^{[2]} \to \PP^{1}$ as unramified covering over
the $4$-punctured $\PP^1$. Denote the loops around 
$0,1,\lambda,\infty \in \PP^1$ by $x_0, x_1, x_2, x_3$ such
that $x_3x_2x_1x_0 =1$. Choose loops around the Weierstra\ss\  
points on $E^{[2]}$ denoted by $c_0,\ldots,c_3$ and  
$a,b \in \pi_1((E^{[2]})^*)$
such that $[a,b]c_3c_2c_1c_0 =1$. This numbering is consistent
with supposing that 
$$h_*(a) = x_3x_2x_1x_3^{-1}, \quad h_*(b)= x_3^2x_2 x_3^{-1}, 
\quad h_*(c_i) = x_i^2.$$
\par
The isogeny $[2]$ doubles the number of unit height cylinders in each 
direction. As $\tilde{\pi}$ is totally ramified over $h^{-1}(\infty)$
all the maximal cylinders of $\pi$ in each direction 
have height $2$.
What we need to ensure is hence that the monodromy images
$m(a), m(ab), m(b) \in S_{d/4}$ 
corresponding to $\tilde{\pi}$ consist of $1$ (resp.\  $d/8$,
resp.\ $r\not\in \{1,d/8\}$) cycles.
\newline
By construction of $\tilde{\pi}$ as fibre product the monodromy $m(c)$ for
$\tilde{\pi}$, where $c \in \pi_1((E^{[2]})^\ast)$ coincides
with the monodromy of $h_*(c)$ for $\beta$.
\newline
Since the monodromy of $x_2$ is trivial, the number of cycles
of $m(a)$ equals the number of preimages $\beta^{-1}(1)$, 
which was $d/8$. 
The number of cycles of $m(b)$ equals the number 
of preimages $\beta^{-1}(\infty)$, which was one. 
Finally, the number of cycles of $m(ab)$ equals the number $r$
of preimages of $0$.
\newline
Going back to the construction of $\beta$, we show that we may
choose $r \not\in \{1,d/8\}$. Suppose $\beta_0$ is a Belyi morphism
as constructed in \cite{Sn94}, totally ramified over $\infty$.
If the number of preimages of $0$ and $1$ does not sum up to
$\deg(\beta)/2$, we take $\beta=4\beta_0(1-\beta_0)$. Otherwise
we take $\beta_1 = x^2 \circ \beta_0$ to ensure this condition
and then take $\beta=4\beta_1(1-\beta_1)$. 
\hfill $\Box$
\par

\begin{Thm} \label{faithfulThm}
For each $\sigma \in G_\QQ$ there is an origami curve $C$ isomorphic
to $\PP^* $, such that $\sigma$
acts non-trivially on the map
$j: C  \to \mrC_{g,[n]}$, 
more precisely such that $j(C) \neq  j^\sigma(C)$.
\end{Thm}
\par
{\bf Proof:}
Take $C$ as constructed in the preceding lemmas and assume
that $j(C) =  j^\sigma(C)$. We claim that this is equivalent to
$j$ being is defined over $K$. To prove this, it is equivalent 
to show that  $$j^\sigma \circ \sigma_{C} = 
\sigma_{\mrC_{g,[n]}} \circ j$$ 
using descent for morphisms between varieties
defined over $K$. The assumption implies that there exists
an automorphism $\varphi$ of $j(C)$, such that 
$$ \varphi \circ j^\sigma \circ \sigma_{C} = 
\sigma_{\mrC_{g,[n]}} \circ j.$$ 
By Lemma \ref{Belyiadjust} and Remark \ref{Monocylrem} 
the stable curves corresponding to the cusps of $C$ have
pairwise distint number of nodes. This implies that
$j(C) \cong \PP^\ast$. Furthermore the number of nodes
of a singular fibre is $G_\QQ$-invariant. Hence $\varphi$
has to be the identity and this proves the claim.
\par 
Take a $K$-rational point $x$ of $C$ and denote the 
corresponding curve by $X$. The above claim
implies the existence of an isomorphism $\sigma_X: X \to X^\sigma$.
By construction
we still have the morphisms $\pi: X \to E$ and $\pi^\sigma: X^\sigma
\to E^\sigma$ plus the canonical morphism $\sigma_E: E \to E^\sigma$.
If we knew that $\pi^\sigma \circ \sigma_X = \sigma_E \circ \pi$
then Lemma \ref{weakfaithful} would lead to a contradiction.
\par
On $X$ (resp.\ $X^\sigma$) we have the differentials 
$\omega_X = \pi^* \omega_E$ (resp.\ $\omega
_{X^\sigma} = \pi^* \omega_{E^\sigma}$).
Our second claim is that $\omega_X = \sigma_X^* \omega_{X^\sigma}$ 
(up to a multiplicative constant).
\newline
The map $\sigma_{\mrC_{g}}$ not only maps $[X] \in (\mrC_{g})_{\ol{\QQ}}$
to $[X^\sigma] \in (\msC_{g})_{\ol{\QQ}}^\sigma$ but also
$C$ to $C^\sigma$. Hence the tangent vector $t_X$ at $[X]$ to $C$
is mapped to the tangent vector $t_{X^\sigma}$ at $[X^\sigma]$ to
$C^{\sigma}$. Now tensor the whole situation by $\CC$ for some
fixed embedding $\ol{\QQ} \to \CC$. The Teichm\"uller metric
induces a (non-linear, but functorial) duality between the 
projectivised tangent and cotangent space of to $\mrC_g$ at $X$.
By construction $(\omega_X)^2$ is the unique (up to scalar multiple)
quadratic differential corresponding to $t_X$ via this duality.
The same is true for $X^{\sigma}$, $t_{X^\sigma}$ and 
$(\omega_{X^\sigma})^2$. 
This implies that $w_X^2 = (\sigma_X^* \omega_{X^\sigma})^2$
and establishes the second claim.
\par 
The periods of $\omega_X$ and $\omega_{X^\sigma}$ define 
a lattice in $\CC$. Integration defines the morphisms
$$f_1: X \to E_1= \CC/\Per(\omega_X), \quad 
f_2: X^\sigma \to E_2= \CC/\Per(\omega_{X^\sigma}).$$ The 
second claim implies that the periods of $\omega_X$ and
$\omega_{X^\sigma}$ are equal. Hence after suitable 
translation there exists
an isomorphism $\sigma_{12}: E_1 \to E_2$ such that
$f_2 \circ \sigma_X = \sigma_{12} \circ f_1$.
\newline
By construction there are isogenies $i_1: E_1 \to E$
and $i_2: E_2 \to E^\sigma$ such that $i_1 \circ f_1 = \pi$
and $i_2 \circ f_2 = \pi^\sigma$. All we need is to
verify that $\sigma_E \circ i_1 = i_2 \circ \sigma_{12}$. 
But we can always lift $\sigma_E$ to an isomorphism $E_1 \to E_2$, 
which has to coincide with $\sigma_{12}$ if $E_1$ (or
equivalently $E$) has no complex multiplication. As elliptic
curves with CM have algebraic integers as $j$-invariants, we
can exclude this by a suitable choice of $x \in C(K)$. 
\hfill $\Box$
\par

\section{A short review of $\widehat{GT}$} \label{onGT}

We start with a summary on braid groups and mapping class groups:
\newline
Denote by $B_n$ the (full) braid group on $n$ strands
and by $P_n$ the pure subgroup, the kernel of the morphism
$p: B_n \to S_n$. Let $\tau_i$ ($i=1,\ldots,n-1$) denote the standard
(Artin) braid generators. We also need the following
elements.
$$ \begin{array}{ll}
y_i = \tau_{i-1}\cdot\ldots\cdot\tau_1\tau_1\cdot\ldots\cdot
\tau_{i-1}, \quad  & i=2,\ldots,n \\ 
w_i = y_2 \ldots y_i, \quad &i=2,\ldots,n.\\
x_{ij} = x_{ji} = (\tau_{j-1}\cdots\tau_{i+1}) \tau_i^2 
(\tau_{j-1}\cdots\tau_{i+1})^{-1}\quad & 1\leq i< j \leq n
\end{array}$$
We denote by $\Gamma_{g,[n]}$ (resp.\ $\Gamma_{g,n}$)
the mapping class group  of a Riemann surface
of genus $g$ with $n$ unordered (resp.\ ordered) points.
For a sphere $\Gamma_{0,[n]}$ equals $B_n/\langle w_n, y_n \rangle$ 
and this presentation
is still valid, when we pass to the profinite completion
$\widehat{\Gamma_{0,n}} = \pi_1(\msC_{0,n})$.
\newline 
We remark that only the {\em orbifold} fundamental
groups of moduli stacks are mapping class groups and this is 
the reason to keep track of the orbifold 
structure of origamis.
\par
In the sequel we use what is called '$\sigma$-convention' in
the appendix of \cite{LoNaSc03}, although we call the standard generators of
the braid groups $\tau_i$ and use $\sigma$ for elements of $G_\QQ$:
paths and braids are composed from the right to the left. Recall
that we abbreviate $\PP^* := \PP^1\sms\{0,1,\infty\}$ and let 
$$x,y,z \in \pi_1^{\topo}(\PP_{\ol{\QQ}}^*, \vec{01}) \subset 
\pi_1(\PP_{\ol{\QQ}}^*,\vec{01}) \cong \widehat{F_2}$$
denote the loops around $0$, $1$ and $\infty$ 
based at the tangential
base point $\vec{01}$ (see below), such that  $xyz = 1$. We denote their
images in the algebraic fundamental group by the same letter.
When using inner automorphisms, the exponent $-1$ is on the
left.
\par
The ($\QQ$-rational) tangential base point $\vec{01}$ defines a splitting
of the exact sequence
$$ 1 \to \pi_1(\PP^*_{\ol{\QQ}}, \vec{01}) \to \pi_1(\PP^*, \vec{01})
\to G_\QQ \to 1$$ 
and with respect to this section the conjugate action of $G_\QQ$
on $\pi_1(\PP^*_{\ol{\QQ}}, \vec{01})$ is 
$$ \begin{array}{lll} x &\mapsto & x^{\chi(\sigma)} \\
y &\mapsto & f_\sigma(x,y)^{-1} y^{\chi(\sigma)} f_\sigma(x,y) \end{array}
\eqno(*),$$
where $\chi(\sigma) \in \widehat{\ZZ}$ is the cyclotomic character and 
$f_\sigma \in \widehat{F_2}$.
Actually $f_\sigma$ lies in the derived subgroup $(\widehat{F_2})'$ 
and $(\chi(\sigma),f_\sigma)$ are known to satisfy the following equations:
$$\begin{array}{rl} (I) & f_\sigma(x,y)f_\sigma(y,x)=1 \\
(II) & f_\sigma(z,x)z^m f_\sigma(y,z)y^m f_\sigma(x,y) x^m =1, 
\,\text{where} \,z = (xy)^{-1}, \,
m=(\chi(\sigma) -1)/2 \\
(III) & f_\sigma(x_{12},x_{23})f_\sigma(x_{34},x_{45})f_\sigma(x_{51},x_{12})
f_\sigma(x_{23},x_{34})f_\sigma(x_{45},x_{51}) =1 \quad \in 
\widehat{\Gamma_{0,5}} \\
\end{array}.$$
We thereby used the convention that $f(a,b)$ denotes the image
of $f$ under the morphism defined by $x \mapsto a, y \mapsto b$.
One defines $\ul{\widehat{GT}}$ as the set of elements
$F= (\lambda,f) \in \widehat{\ZZ} \times \widehat{({F_2})'}$ that
satisfy $(I)$, $(II)$ and $(III)$. $F$ defines an
endomorphism of $\widehat{F_2}$ via $(*)$ and composition
of endomorphisms makes $\ul{\widehat{GT}}$ into a monoid. 
We define $\widehat{GT}$ as the group of invertible elements
of $\ul{\widehat{GT}}$.
\par
Recently several other relations satisfied by the image
of $G_\QQ$ in $\widehat{GT}$ have been found, in particular
in \cite{NaSc00} the following relation in $\widehat{B_3}$
$$\begin{array}{llcl}(IV) &
f_\sigma(\tau_1,\tau_2^2) & = & \tau_2^{4\rho_2(\sigma)}f_\sigma(\tau_1^2,
\tau_2^2)\tau_1^{2\rho_2(\sigma)} (\tau_1\tau_2^2)^{-2\rho_2(\sigma)} \\
&&=&\tau_2^{-4\rho_2(\sigma)}f_\sigma(\tau_1,\tau_2^4)\tau_1^{-2\rho_2(\sigma)}
(\tau_1\tau_2^2)^{2\rho_2(\sigma)}. \\
\end{array}
$$
Here $\rho_p(\sigma)$ denotes the Kummer cocycle on the positive roots
of $\sqrt[n]{p}$ for $n \in \NN$. 
We will make use of the fact that this relation holds in the
subgroup $$\langle \tau_1, \tau_2^2 \mid [\tau_2^2, \tau_1\tau_2^2\tau_1] =1
\rangle \subset \widehat{B_3}.$$ 
\par

\section{Comparison of Galois actions} \label{CompGalAct}

By its very definition we know to express Galois action 
on $\pi_1(\PP_{\ol{\QQ}}^*,\vec{01})$
in terms of $(\lambda,f)$.
On the other hand, from \cite{HaLoSc00} and \cite{NaSc00} 
we 'know' the conjugate action of $G_\QQ$ on $\pi_1(\msC_g)$ with respect 
to some tangential base points based at a maximally degenerate point.
\newline
We can thus use the two-steps origami to compare Galois actions: 
its origami
curve is rational with $3$ cusps, the inclusion $j: C^\orb \to \msC_2$ 
is defined 
over $\QQ$ and the extension $j: \ol{C}^{\rm orb} \to \ol{\msC_{2}}$ 
goes through a maximally degenerate point. Here $\ol{\msC_{2}}$ denotes
the moduli stack of stable curves.
\newline 
We may ignore the orbifold structure of $C^{\orb}$ (see Remark
\ref{orbs2}) for our considerations
because the sequence 
$$ 1 \to (\ZZ/2\ZZ)^2 \to \pi_1(C^\orb, \vec{01}) \to \pi_1(\PP^*,\vec{01})
\to 1 $$
is split. We simply define the $\widehat{GT}$-action on
$\pi_1(C^\orb, \vec{01})$ simply via a fixed splitting.
\par
Consider the loops $a_i$ and $e$ drawn in figure $4$. 
Denote by greek letters (i.e. $\alpha_i,\ve \in \widehat{\Gamma_{2,0}}$) 
the corresponding Dehn twists.
\par
\begin{Thm} \label{newGTrel}
The element $(\lambda,f) \in \widehat{GT}$ respects
the Galois actions on the morphism $j: C^\orb \to \msC_2$ induced from
the two-steps origami ${\mathcal S}_2$ (see figure $4$) 
if and only if
$$ \begin{array}{c}
f(\alpha_3, (\alpha_1^2 \alpha_2)^4)f(\alpha_1^2,\alpha_2^2)
f(\alpha_5^2,\alpha_4^2)\alpha_1^{-4\rho_2(\sigma)}
\alpha_3^{-2\rho_2(\sigma)}\alpha_5^{-4\rho_2(\sigma)} 
f(\alpha_2\alpha_4,\alpha_1^2\alpha_3 \alpha_5^2) \\
= x_{16}^{2\rho_2(\sigma)}
(x_{46}x_{56})^{\rho_2(\sigma)} 
(x_{12}x_{13})^{\rho_2(\sigma)}
(\alpha_2 \alpha_4)^{2\rho_2(\sigma)}
(\alpha_3\alpha_2\alpha_4)^{-2\rho_2(\sigma)}
\end{array} \eqno({\mathcal S}_2)$$
holds in $\widehat{\Gamma_{2,0}}$. The elements $(\lambda,f)$ 
satisfying this relation form a subgroup of $\widehat{GT}$
containing $G_\QQ$.
\end{Thm}
\par
We will prove the theorem in the rest of this section.
\newline
We show now that equation $({\mathcal S}_2)$ is equivalent to the commutativity
of the following diagram:
$$ \xymatrix { \pi_1(C^{\orb}, \vec{01}) \ar[r]^{j_*} 
\ar[d]_{F(C^{{\rm orb}})} &
\pi_1(\msC_{2,0}, \vec{a_*}) \ar[d]^{F(\msC_{2,0})} \\
\pi_1(C^\orb, \vec{01}) \ar[r]^{j_*} &
\pi_1(\msC_{2,0}, \vec{a_*}) \\
}\eqno(*)$$
Here $F=(\lambda, f) \in \widehat{GT}$ and $F(C^\orb)$ (resp.\ $F(\msC_{2,0})$)
are the induced automorphisms on the orbifold fundamental group of $C^\orb$
(resp.\ of $\msC_{2,0}$) as explained
in the previous section (resp.\ as will be explained in section 
\ref{GTonM20}). $\vec{a_*}$ is a base point that will be coveniently
chosen below.
\newline
Once we have shown this, the subgroup property is automatic:
If $F,G \in \widehat{GT}$ satisfy ${\mathcal S}_2$ then
$$j_* \circ (G(C^\orb) \circ F(C^\orb)) = G(\msC_{2,0}) \circ j_*
\circ F(C^\orb) =  G((\msC_{2,0}) \circ F(\msC_{2,0})) \circ j_*$$
and similarly one checks that the subset making the diagram commutative 
is closed under inversion.
\par

\subsection{Comparing tangential base points}
We choose the coordinate $t$ on the origami curve $C \cong \PP^*$ such 
that $0$ 
corresponds to a direction in which the trajectories of the 
Strebel differential $\pi^* \omega$ decompose the surface into
$3$ cylinders. The base point $\vec{a_*}$ of the $G_\QQ$-action on $\msC_2$
will hence correspond to the maximally degenerate point
obtained by shrinking $a_1$, $a_3$ and $a_5$ in figure $4$. 
\par
Somewhat more precisely: A $\pni$-diagram is a trivalent graph
corresponding to a stable curve (maybe with marked points). 
Vertices correspond
to $\PP^1$'s, edges to normal crossings 
and 'loose ends' of the graph correspond to marked points.
To each such diagram one can associate (see \cite{IhNa97}) 
a versal deformation of the stable curve over the
ring $\QQ[[q_1,\ldots,q_{3g-3+n}]]$. 
Taking $q_1=\ldots =q_{3g-3+n}=:q$ gives a $\QQ[[q]]$-valued
point of $\msC_{g,n}$, whose generic fibre we call a {\em standard
tangential base point} associated with the $\pni$-diagram.
It is uniquely determined by the diagram up
to the choice of signs of the $q_i$, or equivalenty what
is called a quilt over the corresponding pants decomposition
in \cite{NaSc00}.
\newline
In our situation there is a unique quilt such that the
resulting tangential base points $\vec{a_*}$ and $j_*(\vec{01})$
are linked by a real path $\gamma$. Thus
the sections $G_\QQ \to \pi_1(\msC_{2},j(\vec{01}))$ induced by
these tangential base points are related by
$$ \gamma^{-1} s_{\vec{a_*}}(\sigma) \gamma
= w_2^{d_l(\sigma)}x_{34}^{d_m(\sigma)} 
w_4^{d_r(\sigma)}   s_{j(\vec{01})}(\sigma) ,$$
where $d_l, d_r$ and $d_m$ are products of Kummer cocycles.
Indeed the morphism $\Sp \QQ[[t]] \to \Sp \QQ[[q_1,q_2,q_3]]$ 
to the base of the
versal deformation induced by $j$ is given by $3$ power series
$p_i(t)$ with leading coefficients say $c_i$ ($i=1,2,3$).
If $q_1$ corresponds to $w_2$ and $c_1 = \prod_p p^{n_p}$, then
$d_l = \sum_p n_p \rho_p(\sigma)$ and similarly for $d_m$ and $d_r$.
\newline
The symmetry of the origami implies that $d_l = d_r$. We
will obtain $d_l=d_r=-2\rho_2(\sigma)$ and $d_m = -\rho_2(\sigma)$  
below automatically by group theory.
\par
We remark that there is a geometric way to determine these 
exponents by comparing tangential base points in 
$\msC_{0,6}$ and $\msC_{2,0}$ using the methods of
\cite{Na02} and \cite{Ic99}. 
\par
 
\subsection{$G_\QQ$-action on $\Gamma_{2,0}$} \label{GTonM20}

We use the notion of $A$- and $S$-moves from \cite{HaLoSc00}.
There is an $A$-move between the pants decompositions of $\Sigma_2$ given
by $\{a_1, a_3, a_5\}$ and $\{a_1, e, a_5\}$. Two $S$-moves change this
into $\{a_2,e,a_5\}$ and into $\{a_2,e,a_4\}$. Hence by \cite{Ic99}
and \cite{NaSc00} the conjugate action of $s_{\vec{a_*}}$ on
$\pi_1(\msC_2,\vec{a_*})$ is given by
$$ \begin{array}{l} \alpha_1 \mapsto \alpha_1^{\chi(\sigma)}, \quad
\alpha_3 \mapsto \alpha_3^{\chi(\sigma)}, \quad
\alpha_5 \mapsto \alpha_5^{\chi(\sigma)} \\
\alpha_2 \mapsto f(\alpha_3, \ve)^{-1} f(\alpha_1^2, \alpha_2^2)^{-1}
\alpha_2^{\chi(\sigma)} f(\alpha_1^2, \alpha_2^2) f(\alpha_3, \ve)\\
\alpha_4 \mapsto f(\alpha_3, \ve)^{-1} f(\alpha_5^2, \alpha_4^2)^{-1}
\alpha_4^{\chi(\sigma)} f(\alpha_5^2, \alpha_4^2) f(\alpha_3, \ve) \\
\end{array}
$$
We should have written $f_\sigma$ instead of $f$, but we will
drop the subscript for simplicity.
Using the topology of the origami, we see that
$j_*$ maps $x,y \in \pi_1(C^\orb,\vec{01})$ to $\alpha_1^2 \alpha_3 
\alpha_5^2$ and
$\alpha_2 \alpha_4$ respectively. Hence the induced 
action on $\alpha_2 \alpha_4 \in \pi_1(\msC_2,
j(\vec{01}))$ is
$$ \alpha_2 \alpha_4 \mapsto f(\alpha_1^2 \alpha_3 \alpha_5^2, \alpha_2 
\alpha_4)^{-1} (\alpha_2 \alpha_4)^{\chi(\sigma)}
f(\alpha_1^2 \alpha_3 \alpha_5^2, \alpha_2 \alpha_4) .$$
Together with the comparison of the tangential base points we obtain
that the diagram $(*)$ commutes, if and only if
$$ f(\alpha_3, \ve) f(\alpha_1^2, \alpha_2^2) f(\alpha_5^2, \alpha_4^2)
w_2^{d_l}x_{34}^{d_m}w_4^{d_r}
f(\alpha_2 \alpha_4, \alpha_1^2 \alpha_3 \alpha_5^2)
\in \Centr_{\widehat{\Gamma_{2,0}}}(\alpha_2\alpha_4).$$
It remains to write this as an equation and to determine $d_l$ and $d_m$.
We may project this expression to $\widehat{\Gamma_{0,[6]}}$, sending
$\alpha_i$ to $\tau_i$ for $i=1,..,5$. More precisely the image lies in 
$p^{-1}(S_4)$, where $p: \widehat{\Gamma_{0,[6]}} \to
S_6$ is induced from $p: B_6 \to S_6$ 
and where $S_4 \subset S_6$ fixes the first and the last marked point. 
We may hence reduce mod $\langle \tau_1^2, \tau_5^2\rangle$ 
(and shift indices by $-1$). Using 
$\ve = w_3^2$ we obtain in $\widehat{\Gamma_{0,[4]}}$
$$ f(\tau_2, \tau_1^4) \,\tau_2^{2d_l} 
f(\tau_1 \tau_3, \tau_2) \in 
\Centr_{\widehat{\Gamma_{0,[4]}}}(\tau_1\tau_3).$$
\par
This is a consequence of relation $(IV)$ in $\widehat{B_3}$. 
Indeed the elements $\tau_1 \tau_3$ and $\tau_2$ satsify 
the defining relation $[\tau_1 \tau_3, \tau_2 \tau_1 \tau_3 \tau_2] =1$.
Hence applying $(IV)$ we obtain in $\widehat{\Gamma_{0,[4]}}$ 
$$ f(\tau_2, \tau_1^4) \,\tau_2^{-2\rho_2(\sigma)} 
f(\tau_1 \tau_3, \tau_2) (\tau_2 \tau_1 \tau_3)^{2\rho_2(\sigma)}
(\tau_1 \tau_3)^{-2\rho_2(\sigma)} =1 \eqno(1)$$
and in particular $d_l = -\rho_2(\sigma)$.
\par

\subsection{Some lemmas on centralizers}

We now prove some lemmas that will be applied in the next
section. The author is grateful to W.~Herfort for 
these results.
\par  
\begin{Lemma} \label{freeonbasis}
Let the profinite group $G$ act continuously and freely on a profinite
space $X$. Then for the induced action on $\free X$ the equation
$\varphi f=f$ for some $\varphi \in G$ and $f\in \free X$ 
yields either $\varphi=1$ 
or $f=1$.
\end{Lemma} 
\par
{\bf Proof:} 
Suppose, on the contrary,  
there exist non-trivial elements $\varphi \in G$ and $f\in \free X$
with $\varphi f=f$. 
We first want to show that $G$ and $X$ can be assumed
to be finite.
Since $G$ acts freely on $X$, using
Lemma 5.6.5 (a) in \cite{RiZa00} we find a continuous section 
$\sigma:G\backslash X\rightarrow X$. Hence we may
identify $X$ as a $G$-space with $G\times\Lambda$, where, for short,
we have put $\Lambda:=G\backslash X$, with the left regular action of $G$.
Let $N$ be a normal open subgroup of $G$ and $R$ a clopen relation on 
$\Lambda$,
such that $g\not\in N$ and $f$ has a non-trivial image under canonical
projection from $\free{G\times\Lambda}$ onto $\free{G/N\times\Lambda/R}$.
Such $N$ and $R$ exist as Proposition 1.7 in \cite{GiLi72} shows.
Since $G/N$ acts freely on $G/N\times\Lambda/R$ we have shown that indeed it
suffices to assume $G$ and $X=G\times\Lambda$ both to be finite. 
\newline
Let $\Gamma:=G*\free \Lambda$ be the free product and
define a map $\eta$ from $G\cup\free\Lambda$ 
to the holomorph $H:=G\ltimes\free X$ as follows.
We send
$\varphi \in G$ to $\varphi \in G$ (as a subgroup of the holomorph) and
we extend the map that sends $\lambda\in\Lambda$ to 
$(1,\lambda)\in G\times\Lambda$, to a continuous homomorphism
from all of $\free\Lambda$ to $H$ by using the universality of the free
group $\free\Lambda$.
Use the universal property of $\Gamma$ being a free product in order to
extend $\eta:G\cup\free\Lambda\rightarrow H$ to a continuous epimorphism
$\omega:\Gamma\rightarrow H$. 
\par
We claim $\omega$ to be an isomorphism.
Indeed, since $\omega$ induces the identity on $\Gamma/(\free\Lambda)_\Gamma$,
we conclude $\Ker\;\omega\le(\free\Lambda)_\Gamma$.
(For a profinite group $G$ and a subset $A$ of $G$ let $(A)_G$ denote
the normal closure, i.e., the smallest closed normal subgroup of $G$
containing $A$). Use the Kurosh Subgroup Theorem (e.g.\
Thm.\ 9.1.9 in \cite{RiZa00}) by applying it to the
normal open subgroup $(\free\Lambda)_\Gamma$ 
in order to see that its rank equals $|G|\times |\Lambda|$.
Since $(\free\Lambda)_\Gamma$ is hopfian and it goes onto $\free X$
(viewed as a subgroup of $H$) conclude $\Ker\; \omega=\{1\}$. 
\newline
Let us point out that the action of $\varphi$ on $X$ becomes
conjugation in the holomorph $H$. Now use Theorem 9.1.12 in 
\cite{RiZa00} in order to see that $C_\Gamma(g)\le G$,
so that applying $\omega$ one finds $C_H(g)\le G$. 
This however contradicts the choice of $\varphi$ and $f$.
\hfill $\Box$
\par
We apply this lemma in the following two cases:
\par
\begin{Lemma} \label{F3cent}
Let $\widehat{F_3} = \langle x,y,z\rangle$ be
the profinite free group on three generators and $\varphi$
the following automorphism of $\widehat{F_3}$: $\varphi(x) = x$,
$\varphi(y) = x y x^{-1}$ and $\varphi(z) = x^{-1} z x$. 
Then the fixed group of $\varphi$
in $\widehat{F_3}$ is the profinite free group generated by $x$.
\end{Lemma}
\par
{\bf Proof:}
Let $N$ denote the normal subgroup generated by $y$ and $z$. By
specialising Thm.~8.1.3 in \cite{RiZa00} to the finite case, we 
deduce that $N$ is the profinite free group on the generators
$X= \{ x^{-l}yx^l, x^{-l}zx^l \; \mid \; l \in \widehat{\ZZ}\}$.
$\varphi$ acts freely on $X$ and Lemma \ref{freeonbasis} yields
$\Fix(\varphi) \cap N = \{1\}$. Hence the isomorphism
$\widehat{F_3}/N \to \langle x \rangle$ induces an
isomorphism $\Fix(\varphi) \to \langle x \rangle$. \hfill $\Box$ 
\par
\begin{Lemma} \label{F4cent}
Let $\widehat{F_4} = \langle w,x,y,z\rangle$ be the profinite 
free group on $4$ generators and $\varphi$ the following
automorphisms of $\widehat{F_4}$: $\varphi(w)=w$, $\varphi(x) = x$, 
$\varphi(y) = x y x^{-1}$ and $\varphi(z)= w^{-1} x^{-1} z x w$. 
Then the fixed group of $\varphi$
in $\widehat{F_4}$ is the profinite free group generated by 
$w$ and $x$.
\end{Lemma}
\par
{\bf Proof:} Let $N$ denote the normal subgroup generated by
$z$. As above, $N$ is free on the generators
$X = \{ u^{-1}zu  \mid \; u \in \langle w,x,y \rangle \}$.
We check that $\varphi$ acts freely on $X$.
\newline
Suppose for $l \in \widehat{\ZZ}$ we have 
$\varphi^l(u^{-1}zu) = u^{-1}zu$. By definition of $\varphi$
this implies 
$$(xw)^l\varphi^l(u)u^{-1} \in \Centr_{\widehat{F_4}}(z).$$
This centralizer equals $\langle z \rangle$ and has trivial
intersection with $\langle w,x,y \rangle$. Thus
$(xw)^l\varphi^l(u)=u$. If we consider this modulo the
normal subgroup $N_y$ generated by $y$, the action of $\varphi$
is trivial, hence $(xw)^l \in N_y$. This is only possible for $l=0$. 
\newline
By Lemma \ref{freeonbasis} hence $\Fix(\varphi) \cap N = \{ 1\}$.
Thus $\Fix(\varphi)$ injects into $\widehat{F_4}/N \cong \widehat{F_3}
= \langle w,x,y \rangle$. Let $N_y \subset \widehat{F_3}$
be the normal subgroup generated by $y$. With the same
arguments we conclude that $\Fix(\varphi) \cap N_y = \emptyset$
and hence $\Fix(\varphi)=\langle w,x \rangle$.
\hfill $\Box$ 

\par

\subsection{Lifting the relation}

We now want to lift the equation  $(1)$ 
successively to $\widehat{\Gamma_{0,[5]}}$,
$\widehat{\Gamma_{0,[6]}}$ and to $\widehat{\Gamma_{2,0}}$.
Let $S_4 \subset S_5$ be the permutation group of the last $4$ strings
and $p: \widehat{\Gamma_{0,[5]}} \to S_5$ the permutation representation.
Reducing mod $\langle \tau_1^2 \rangle$ gives a morphism
$p_1: p^{-1}(S_4)  \to 
\widehat{\Gamma_{0,[4]}}$. Undoing the shift of indices, we know that
$$ 
\begin{array}{c}
f(\tau_3, (\tau_1^2 \tau_2)^4)f(\tau_1^2,\tau_2^2)
\tau_1^{2d_l} \tau_3^{-2\rho_2(\sigma)} f(\tau_2\tau_4,\tau_1^2\tau_3) 
(\tau_3 \tau_2 \tau_4)^{2\rho_2(\sigma)}(\tau_2 \tau_4)^{-2\rho_2(\sigma)}
\\
\in \Ker(p_1) \cap \Centr_{\widehat{\Gamma_{0,[5]}}}(\tau_2 \tau_4).
\end{array} \eqno(2)$$
\newline
$\Ker(p_1)$ is the free profinite group on the three generators 
$x_{12}, x_{12}x_{13}$ and $x_{15}$ and the conjugate action
of the square of $\tau_2\tau_4$ on $\Ker(p_1)$ is given by
$$\begin{array}{lcl}
(\tau_2\tau_4)^{-2} x_{12}x_{13} (\tau_2\tau_4)^2 &=& x_{12}x_{13} \\
(\tau_2\tau_4)^{-2} x_{12} (\tau_2\tau_4)^2 &=& (x_{12}x_{13})\;\;\; x_{12}\; 
(x_{12}x_{13})^{-1} \\
(\tau_2\tau_4)^{-2} x_{15} (\tau_2\tau_4)^2 &=& (x_{12}x_{13})^{-1}
x_{15}\; (x_{12}x_{13}) \\
\end{array}$$
\par
Using Lemma \ref{F3cent} we conclude that the above expression $(2)$
equals a power of $x_{12}x_{13}$. To determine the exponent, we
cannot simply abelianize the subgroup $p^{-1}((24)(35))$ 
of $\widehat{\Gamma_{0,5}}$ because $x_{12}x_{13}$
vanishes. Therefore we first shift indices by $-1$
and use the natural embedding $\widehat{B_4}/\langle w_4\rangle
\to \widehat{\Gamma_{0,[5]}}$, $\tau_i \mapsto \tau_i$ for
$1\leq i \leq 4$. The equation then becomes
$$ f(\tau_2, x_{34}^2)f(z_3,\tau_1^2)z_3^{d_l}x_{23}^{-\rho_2(\sigma)}
f(\tau_1 \tau_3, z_3 \tau_2)(\tau_2\tau_1\tau_3)^{2\rho_2(\sigma)}
(\tau_1\tau_3)^{-2\rho_2(\sigma)} 
= (z_3 \tau_1 z_3 \tau_1^{-1})^a \eqno(3)$$
where $z_3 = (\tau_2 \tau_3)^3$. Note that in $\widehat{\Gamma_{0,5}}$
we have $x_{12}=(\tau_3\tau_4)^3$. 
We now let $k_5 := [\tau_1 \tau_3, \tau_2 z_3 
\tau_1 \tau_3 \tau_2 z_3]$ and 
map equation $(2)$ 
to $\widehat{G}= \widehat{B_4}/\langle w_4, k_5\rangle$ in order
to apply relation $(IV)$ to $f(\tau_1 \tau_3, z_3\tau_2 )$.
\par
The abelianization of $\widehat{G}$ is generated by
the $x_{ij}$ for $1 \leq i < j \leq 4$ with the relations 
$$x_{12} \equiv (x_{13}x_{14}x_{23}x_{24}x_{34})^{-1} \quad \text{and}
\quad 
x_{13} \equiv  x_{24}$$ due to the factorization mod $w_4$ and $k_5$.
\newline
A simple calculation in $\widehat{G}^{\rm ab}$ yields
$[z_3 \tau_2, (\tau_1 \tau_3)^2] \equiv 
x_{12}^{-1} x_{13}x_{24}  x_{34}^{-1}$.
On this commutator $z_3 \tau_2$ acts by $(-1)$ and $(\tau_1 \tau_3)^2$
acts trivially. Hence Ihara's Blanchfield-Lyndon calculus (see
\cite{NaSc00} Lemma 2.3) implies
$$ f(z_3 \tau_2, (\tau_1 \tau_3)^2) \equiv  
(x_{12}^{-1}  x_{13} x_{24}  x_{34}^{-1})^{-\rho_2(\sigma)}.$$
For the same reason
$$ f(\tau_2, x_{34}^2) \equiv (x_{24}^{2}x_{34}^{-2})^{-\rho_2(\sigma)}.$$
In $\widehat{G}^{\rm ab}$ the right hand side of equation $(3)$ equals
$$ (z_3 \tau_1 z_3 \tau_1^{-1})^a \equiv (x_{14}x_{23}x_{24}^2 x_{34}^2)^{a}.$$
The left hand side, using relation $(IV)$ and the above formula
(note that $f(z_3,\tau_1^2)$ vanishes), equals
$$ LHS(2) \equiv x_{14}^{\rho_2(\sigma)}x_{23}^{d_l+3\rho_2(\sigma)}
(x_{24}x_{34})^{d_l+4\rho_2(\sigma)}.$$
We conclude that $d_l = -2\rho_2(\sigma)$ and $a=\rho_2(\sigma)$.
\par
To go on to $\widehat{\Gamma_{0,[6]}}$ let $p^{-1}(S_4) \subset
\widehat{\Gamma_{0,[6]}}$ be the subgroup that fixes the first and last
marked point. Denote by $p_6:p^{-1}(S_4) \to 
\widehat{\Gamma_{0,[5]}}$ the
reduction mod $\langle \tau_5^2 \rangle$.
\newline
If we undo the indexshift and multiply the expression by 
$(x_{46}x_{56})^{-\rho_2(\sigma)}$ to obtain an expression
symmetric with respect to $\tau_i \leftrightarrow \tau_{6-i}$
(note that $x_{12}x_{13}\leftrightarrow x_{46}x_{56}$)
we conclude 
$$
\begin{array}{l} f(\tau_3, (\tau_1^2 \tau_2)^4)f(\tau_1^2,\tau_2^2)
f(\tau_5^2,\tau_4^2)\tau_1^{-4\rho_2(\sigma)}
\tau_3^{-2\rho_2(\sigma)}\tau_5^{-4\rho_2(\sigma)} 
f(\tau_2\tau_4,\tau_1^2\tau_3 \tau_5^2) \cdot \\
\cdot 
(\tau_3\tau_2\tau_4)^{2\rho_2(\sigma)}
(\tau_2 \tau_4)^{-2\rho_2(\sigma)}(x_{12}x_{13})^{-\rho_2(\sigma)} 
(x_{46}x_{56})^{-\rho_2(\sigma)}
\in \Ker(p_6) \cap \Centr_{\widehat{\Gamma_{0,[6]}}}(\tau_2 \tau_4)
\end{array}. \eqno(4)$$
\par
$\Ker(p_6)$ is the free profinite group on 4 generators
$x_{16}$, $x_{46}x_{56}$, $x_{56}$ and $x_{26}$ 
and the conjugate action of $(\tau_2\tau_4)^2$ (we only
use the action of the square) on $\Ker(p_6)$ is given by
$$\begin{array}{lcl}
(\tau_2\tau_4)^{-2} x_{16} \;(\tau_2\tau_4)^2 &=& x_{16} \\
(\tau_2\tau_4)^{-2} x_{46}x_{56} \;(\tau_2\tau_4)^2 &=& x_{46}x_{56} \\
(\tau_2\tau_4)^{-2} x_{56} \;(\tau_2\tau_4)^2 &=& \phantom{x_{16}^{-1}}
(x_{46}x_{56})\;\;\; x_{56}\; 
(x_{46}x_{56})^{-1} \\
(\tau_2\tau_4)^{-2} x_{26} \;(\tau_2\tau_4)^2 &=& x_{16}^{-1}
(x_{46}x_{56})^{-1}
x_{26}\; (x_{46}x_{56}) x_{16} \\
\end{array}$$
\par
By Lemma \ref{F4cent} and the 
symmetry with respect to $\tau_i \leftrightarrow \tau_{6-i}$ 
the expression $(4)$ equals a power of $ x_{16}$. To determine the
exponent, we use the same technique as above and apply relation
$(IV)$ twice. In order to be able to do so, we consider the
equation in $\widehat{G}=\widehat{\Gamma_{0,6}}/\langle k_6 \rangle$,
where $k_6 = [\tau_2 \tau_4, \tau_1^2 \tau_3 \tau_5 \cdot \tau_2 \tau_4 \cdot
\tau_1^2 \tau_3 \tau_5 ]$. The abelianization of $\widehat{\Gamma_{0,6}}$
is generated by $x_{12},x_{13},x_{14},x_{15},x_{23},x_{24},x_{25},
x_{34},x_{35}$. In $\widehat{G}^{{\rm ab}}$ we have the supplementary
relation $x_{35}^2 \equiv (x_{12}^{-1}x_{13}x_{14}^{-1}x_{15}
x_{24}^{-1})^{-2}$. We use the Blanchfield-Lyndon calculus again
to obtain
$$ f(\tau_2^2\tau_4^2,\tau_1^2\tau_3\tau_5^2) \equiv (x_{12}x_{13}x_{14}x_{24}^2
x_{34}x_{15}x_{25}x_{35}^2)^{\rho_2(\sigma)}$$
and
$$ f(\tau_3,(\tau_1^2 \tau_2)^4) \equiv (x_{12}x_{14}^{-1}x_{24}^{-1}
x_{23})^{2\rho_2(\sigma)}.$$
In $\widehat{G}^{{\rm ab}}$ we can apply relation $(IV)$ to
$f(\tau_2\tau_4,\tau_1^2\tau_3 \tau_5^2)$ and use to
above formula. By direct calculation one verifies that
$$ (\tau_2 \tau_4 \tau_1^2\tau_3\tau_5^2)^{-2\rho_2(\sigma)}
(\tau_3\tau_2\tau_4)^{2\rho_2(\sigma)} \equiv (x_{13}^2 x_{15}^2 
x_{24}^{-1} x_{35})^{-\rho_2(\sigma)}.$$
Hence the left hand side
of equation $(4)$ sums up to
$$ LHS(3) \equiv (x_{12} x_{13}x_{14}x_{15})^{-2\rho_2(\sigma)} 
\equiv x_{16}^{2\rho_2(\sigma)} = RHS(3).$$
\par
Finally we want to lift the equation to $\widehat{\Gamma_{2,0}}$,
which is an extension
of $\widehat{\Gamma_{0,[6]}}$ by the hyperelliptic involution
generated by the central element $w_{5}$. Note that equation $(4)$ 
involves only even powers of $\alpha_1$ and $\alpha_5$. We can reduce 
modulo $\langle \alpha_1^2, \alpha_5^2 \rangle$ 
(and shift indices by $(-1)$) 
to obtain in the sphere braid group $\widehat{H_4}= 
\widehat{B_4}/\langle y_4 \rangle$
$$ f(\tau_2, \tau_1^4) \,\tau_2^{-2\rho_2(\sigma)} 
f(\tau_1 \tau_3, \tau_2) (\tau_2 \tau_1 \tau_3)^{2\rho_2(\sigma)}
(\tau_1 \tau_3)^{-2\rho_2(\sigma)} \in {\rm center(\widehat{H_4})}.$$
We had noticed above that this expression equals $1$ in 
$\widehat{\Gamma_{0,[4]}}$ by using relation $(IV)$. Denote by 
$w_3 = (\tau_1 \tau_2)^3$ the center of $\widehat{H_4}$. Relation $(IV)$
now tells us that 
$$ 
\begin{array}{l}f(\tau_2, \tau_1^4) \,\tau_2^{-2\rho_2(\sigma)} 
f(\tau_1 \tau_3, \tau_2) (\tau_2 \tau_1 \tau_3)^{2\rho_2(\sigma)}
(\tau_1 \tau_3)^{-2\rho_2(\sigma)} = \\
= f(\tau_2, \tau_1^4) f(\tau_1^2 \tau_3^2, \tau_2)
= f(\tau_2, \tau_1^4) f(\tau_1^4 w_3, \tau_2) = 1.
\end{array} $$ 
This completes the proof of the theorem. \hfill $\Box$
\par

Martin M{\"o}ller: Universit{\"a}t Essen, FB 6 (Mathematik) \newline 
45117 Essen, Germany \newline
e-mail: martin.moeller@uni-essen.de \newline

\end{document}